\documentclass[a4paper,10pt]{article}
\usepackage[cp1251]{inputenc}
\usepackage[russian, ukrainian, english]{babel}
\usepackage{amsmath, amsthm, amsfonts, amssymb}

\newtheorem{thm}{Theorem}[section]

\theoremstyle{remark}
\newtheorem{rem}{Remark}[section]
\theoremstyle{definition}
\newtheorem{dfn}{Definition}[section]
\newtheorem{exm}{Example}[section]
\usepackage{comment}
\begin{document}

\begin{center}
{\bf CLARK  REPRESENTATION FORMULA FOR THE SOLUTION TO
EQUATION WITH INTERACTION}
\footnote{
200 {\it Mathematics Subject Classification}: 60H35; 	60H07; 93E03; 393E10.\newline
\qquad{\it Key words and phrases}. Stochastic differential equations with interaction, Clark representation, Clark-Okone formula.
}\\[20pt ]
{\it Jasmina \DJ or\dj evi\'c and Andrey Dorogovtsev}\\[20pt]

{\small {Dedicatory: Jasmina \DJ or\dj evi\'c is supported by STORM-Stochastics for Time-Space Risk Models,\\ granted by Research Council of Norway -- Independent projects: ToppForsk. Project nr. 274410,\\
The Faculty of Mathematics and Natural Sciences, University of Oslo, Norway}.}

\end{center}

{\bf Abstract}.
In this paper Clark-Ocone representation for solution to measure-valued equation with interaction is studied. It is proven that the integrand is absolutely continuous with respect to Lebesgue measure.
\vskip20pt

\section{Introduction}

In this article we consider the measure-valued process, which is the solution to the following one-dimensional equation with interaction
\begin{eqnarray}
&&d x(u, t)=\int_{\mathbb{R}} {a}(x(u, t), x(v, t)) \mu_{0}(d v) d t+{b}(x(u, t)) d W(t)\nonumber\\
&&x(u, 0)=u, \quad u \in \mathbb{R}.\label{1}
\end{eqnarray}
Here $\mu_{0}$ is a probability measure on $\mathbb{R},$ which can be treated as an initial mass distribution of the infinite system of particles whose trajectories are $x(u, t), t \geq 0, u \in \mathbb{R}$.
In such interpretation the measure
$$
\mu_{t}=\mu_{0} \circ x^{-1}(u, t), \quad t \geq 0
$$
is the distribution of the mass of particles at the moment $t \geq 0 .$ Hence, the equation (\ref{1})
 is the partial case of more general equation with interaction
\begin{eqnarray}
&&d x(u, t)={a}\left(x(u, t), \mu_{t}\right) d t+{b}\left(x(u, t), \mu_{t}\right) d W(t)\nonumber \\
&&x(u, 0)=u, u \in \mathbb{R},\nonumber \\
&&\mu_{t}=\mu_{0} \circ x^{-1}(u, t), t \geq 0.\label{2}
\end{eqnarray}
This type of equation was introduced and studied by A. Dorogovtsev in \cite{aa}.

 Here we consider the measure valued process $\mu_{t}$ as a functional of the noise $W(\cdot)$. It is natural question for one to ask, what would be the Clark representation for the random measure $\mu_{t} .$ Let us recall on Clarck representation for random variables.

\begin{thm} (Clark \cite{c}) Let $\alpha$ be square-integrable random variable measurable with respect to a Wiener process $\{W(t) ; t \in[0 ; 1]\} .$ Then there exists non-anticipating random function $\{\xi(t) ; t \in[0 ; 1]\}$ such that
\begin{equation}
\alpha=E \alpha+\int_{0}^{1} \xi(t) d W(t).\label{3}
\end{equation}
The random function $\xi$ is unique up to the stochastic equivalence.

\end{thm}

The representation from previous theorem can be transformed into representation of martingales adapted to the Wiener filtration and it has a wide application in construction of backward stochastic differential equations (see  \cite{moi2}, \cite{moi1}). Thus it would be interesting to find  Clark representation  for measure valued process  $\mu_t$. The idea will be developed in following way. For any test function $\varphi \in C_{0}^{\infty}\left(\mathbb{R}^{d}\right),$  consider integral
$$
\left\langle\varphi, \mu_{t}\right\rangle:=\int_{\mathbb{R}^{d}} \varphi(v) \mu_{t}(d v)=\int_{\mathbb{R}^{d}} \varphi(x(u, t)) \mu_{0}(d u).
$$
This random variable has Clark representation
\begin{equation}
\left\langle\varphi, \mu_{t}\right\rangle=E\left\langle\varphi, \mu_{t}\right\rangle+\int_{0}^{t} \Theta_{\varphi}(s) d W(s).\label{4}
\end{equation}
It is obvious that $\Theta_{\varphi}$ is linear with respect to $\varphi$. The aim is to explore if $\Theta_{\varphi}$ can be expressed as
$$
\Theta_{\varphi}(s)=\left\langle\varphi, \varkappa_{s}^{\varphi}\right\rangle,
$$
where $\varkappa_{s}^{\varphi}$ is a random signed measure. So we are looking for the Clark representation formula of the following kind

\begin{equation}
\mu_{t}=E \mu_{t}+\int_{0}^{t} \varkappa_{s}^{\varphi} d W(s),\label{5}
\end{equation}
which is valid in the sense of action on the test function.


In order to deduce representation (\ref{5}), we will use Clark-Ocone's formula (see \cite{b}). This formula provides the expression for integrand in Clark representation in terms of stochastic derivative for the random variable $\alpha$.

Namely, if $\alpha$ has stochastic derivative $D \alpha,$ then
$$
\alpha=E \alpha+\int_{0}^{t} E\left(D \alpha(s) \mid \mathcal{F}_{s}\right) d W(s),
$$
where $\mathcal{F}_{s}=\sigma\{W(l), l \leq s\}.$

This is a part of general question. Let $\mathcal{B}$ be a real separable Banach space, and let $\xi \in \mathcal{B}$ be a random element with a finite second moment of the norm. Suppose that $\xi$ is measurable with respect to a Wiener process $W(t), t \in[0,1]$. Then
for arbitrary $\varphi \in \mathcal{B}^{*}\left(\mathcal{B}^{*}\right.$ is dual space of $\left.\mathcal{B}\right)$ random variable $\langle\varphi, \xi\rangle$ (scalar product) has Clark representation
$$
\langle\varphi, \xi\rangle=\langle\varphi, \bar{\xi}\rangle+\int_{0}^{1} R_{\varphi}(t) d W(t),
$$
where $\bar{\xi}$ is the (Bochner) mean value of $\xi$ and $R_{\varphi}(\cdot)$ is the square-integrable integrand which linearly depends on $\xi$. The main question can be formulated in the following way: $\operatorname{Can} R_{\varphi}$ be represented as
$$
R_{\varphi}=\langle\varphi, \eta(t)\rangle,
$$
where $\eta(t), t \in[0,1]$ is $\mathcal{B}$ -valued random process, which is adapted to the Wiener filtration. Answer to this question depends on the geometry of the space $\mathcal{B}$.

\begin{dfn} (See \cite{ak}.) Let $\mathcal{B}$ be a Banach space. It has unconditional martingale (shorter UMG) property if for every martingale difference $d_{1}, d_{2}, \ldots, d_{n}$ with the finite $p$ -th moment of the norm (for $p>1$ ) and arbitrary choice of signs $\varepsilon_{k}=\pm 1,$ next inequality holds
$$
E\left\|\sum_{k=1}^{n} \varepsilon_{k} d_{k}\right\|^{p} \leq C_{p} E\left\|\sum_{k=1}^{n} d_{k}^{p}\right\|
$$
where constant $C_{p}$  depends only on $p$.
\end{dfn}

It was prove \cite{ak} that if the property from Definition 1.1 holds for $p>1$, then it holds dor all $p>1$. Per example, spaces $L^{p}$, for $p>1$ with respect to $\sigma$-finite measure have the UMG property, while $L^{1}$ does not have this property. It occurs that UMG property has a crucial role for Clark representation. The answer to introduced problem is positive for Banach spaces which have UMG property. As it can be easily seen the space of finite signed measures on $\mathcal{R}$ with distance in variation contains $L^{1}$ as a subspace. But in certain cases Clark representation for random measures still  have good properties. Let us consider following example.

\begin{exm}
\label{exm1.1}
 Let $W(t), t \in[0,1]$ be a Wiener process. Consider a random measure $\delta_{W(1)} .$ Then for function $f \in C_{0}^{\infty}(\mathbb{R})$ one can get
$$
\begin{aligned}
\left\langle f, \delta_{W(1)}\right\rangle &=f(W(1))=\int_{\mathbb{R}} f(u) p_{1}(u) d u+\int_{0}^{1} E\left(D f(W(1))(t) \mid \mathcal{F}_{t}^{W}\right) d W(t) \\
&=\int_{\mathbb{R}} f(u) p_{1}(u) d u+\int_{0}^{1} E\left(f^{\prime}(W(1)) \mid \mathcal{F}_{t}^{W}\right) d W(t)
\end{aligned}
$$
$$
=\int_{\mathbb{R}} f(u) p_{1}(u) d u+\int_{0}^{1} \int_{\mathbb{R}} f'(W(t)+v) p_{1-t}(v) d v d W(t)
$$
$$
=\int_{\mathbb{R}} f(u) p_{1}(u) d u+\int_{0}^{1} \int_{\mathbb{R}} f(W(t)+v)\left(-p_{1-t}^{\prime}(v)\right) d v d W(t).
$$

Hence, in the sense of pairing with the test functions, the following formula holds
$$
\delta_{W(1)}=p_{1}-\int_{0}^{1} p_{1-t}^{\prime}(-W(t)) d W(t).
$$
\end{exm}

\begin{rem} Stochastic derivative for solution of equation with interaction was discussed in \cite{kh}, where also the probability density was established. For our purpose we need the estimation of the density  which will lead to completely different approach than the one which is used in \cite{kh}.
\end{rem}

\section{ Stochastic derivative of the solution }
As it was already mentioned in the Introduction, stochastic derivative of  $\alpha$ will be denoted with $D \alpha .$ Further, let us suppose that the coefficients ${a}(\cdot, \cdot), {b}(\cdot),$ of equation (\ref{1}) have continuous bounded derivatives up to second order respect to all variables. Then following statement holds.

\begin{thm} The solution $x(u, t)$ of the equation (1) has the stochastic derivative $\eta_s(u, t):=D x(u, t)(s),$ which satisfies equation

\begin{eqnarray}
&&d \eta_{s}(u, t)= \int_{\mathbb{R}} {a}(x(u, t), x(v, t)) \mu_{0}(d v) \eta_{s}(u, t) d t \nonumber \\
&&\phantom{d \eta_{s}(u, t)=}+\int_{\mathbb{R}} {a}(x(u, t), x(v, t)) \mu_{0}(d v) \eta_{s}(v, t) d t \nonumber \\
&&\phantom{d \eta_{s}(u, t)=}+{b}(x(u, t)) \eta_{s}(u, t) d W(s)\nonumber \\
&&\eta_{s}(u, s)=b(x(u, s)).\label{6}
\end{eqnarray}

Moreover, for every $T>0$ and $p>0$
$$
\sup _{s \leq t \leq T, u \in \mathbb{R}} E\left|\eta_{s}(u, t)\right|^{p}<+\infty.
$$
\end{thm}
\begin{proof}  We will recall on the iteration sequence for the solution of eq.(\ref{1}). Let
$$
\begin{array}{l}
x^{0}(u, t)=u, u \in \mathbb{R}, t \geq 0 \\
\mu_{t}^{0}=\mu_{0} \circ x^{0}(\cdot, t)^{-1}=\mu_{0}.
\end{array}
$$
For $n>0$ put
$$
d x^{n+1}(u, t)=\int_{\mathbb{R}} {a}\left(x^{n}(u, t), x^{n}(v, t)\right) \mu_{0}(d v) d t+{b}\left(x^{n}(u, t)\right) d W(t)
$$
$$
x^{n+1}(u, 0)=u, \quad u \in \mathbb{R} .
$$
It has been proved in \cite{aa}, that this sequence converges to a solution $x(\cdot, \cdot)$ of eq. (\ref{1}), in a sense that for every $K, T, p>0$
$$
E \sup _{|u| \leq K, 0 \leq t \leq T} E\left|x^{n}(u, t)-x(u, t)\right|^{p} \longrightarrow 0, \quad n \longrightarrow+\infty .
$$
Using this, existence of stochastic derivative of process $x(\cdot, \cdot)$ can be proved.

Denote
$$
\eta_{s}^{n}(u, t)=D x^{n}(u, t)(s).
$$

It can be checked by induction that $\eta_{s}^{n}(u, t)$ exists and  satisfies the relation
$$
d \eta_{s}^{n+1}(u, t)=\int_{\mathbb{R}} {a}_{1}'\left(x^{n}(u, t), x^{n}(v, t)\right)  \eta_{s}^{n}(u, t) \mu_{0}(d v)d t
$$
$$
\begin{aligned}
&+\int_{\mathbb{R}} {a}_{2}'\left(x^{n}(u, t), x^{n}(v, t)\right)  \eta_{s}^{n}(v, t) \mu_{0}(d v)d t \\
&+{b}'\left(x^{n}(u, t)\right) \eta_{s}(u, t) d W(s), \\
\eta_{s}^{n+1}(u, s)=b &\left(x^{n}(u, s)\right).
\end{aligned}
$$
We have used the know result that for the continuous, stochastically differentiable adapted process 
$\beta(r), r \geq 0,$ the stochastic derivative of It\^o integral
$$
D \int_{0}^{t} \beta(r) d W(r)(s)=\beta(s)+\int_{s}^{t} D \beta(r)(s) d W(r).
$$
The proof of the theorem follows from the fact that the stochastic derivative is a closed operator.
The estimation of the moments of stochastic derivative follows from the Gronwall Bellman lemma.

 \end{proof}

\begin{rem}
\label{rem2.1} Repeating the same steps, it can be shown that a second stochastic derivative of process $x(\cdot, \cdot)$ exists and that it has bounded moments.
\end{rem}

\begin{rem}
\label{rem2.2}
It should be noted that in \cite{ak} author provided the estimates for stochastic derivative under weaker conditions.
\end{rem}

Note that in case  $b(x)\equiv 1, x\in \mathbb{R} $  the stochastic derivative satisfies the equation without the stochastic differential and can be bounded from zero and infinity.

 So one can try to apply Clark-Okone formula to $\mu_{1}$.

$$
\begin{aligned}
\left\langle\varphi, \mu_{1}\right\rangle &=\left\langle\varphi, E \mu_{1}\right\rangle+\int_{0}^{1} E\left(D\left\langle\varphi, \mu_{1}\right\rangle(t) \mid \mathcal{F}_{t}^{W}\right) d W(t) \\
&=\left\langle\varphi, \bar{\mu}_{1}\right\rangle+\int_{0}^{1} E\left(D\left\langle\varphi(x(\cdot, 1)), \mu_{0}\right\rangle(t) \mid \mathcal{F}_{t}^{W}\right) d W(t) \\
&=\left\langle\varphi, \overline{\mu_{1}}\right\rangle+\int_{0}^{1} \int_{\mathbb{R}} E\left(\varphi^{\prime}(x(u, 1)) \eta_{t}(u, t) \mid \mathcal{F}_{t}^{W}\right) d \mu_{0}(u) d W(t).
\end{aligned}
$$

Let us consider expression under the integrals, with the evident notations

\begin{eqnarray*}
&&E\left(\varphi^{\prime}(x(u, 1)) \eta_{t}(u, t) \mid \mathcal{F}_{t}^{W}\right) \\
&&\quad=E\left(\varphi^{\prime}(x(u, 1))\left(1+\int_{t}^{1} \alpha(u, s) \eta_{t}(u, s) d s+\int_{t}^{1} \int_{\mathbb{R}} \beta(u, v, s) \eta_{t}(v, s) \mu_{0}(d v) d s\right) \mid \mathcal{F}_{t}^{W}\right).
\end{eqnarray*}

Now,
$$
E\left(\varphi^{\prime}(x(u, 1)) \mid \mathcal{F}_{t}^{W}\right)=E\left(\varphi^{\prime}\left(x_{t, 1}\left(\mu_{t}, x(u, t)\right)\right) \mid \mathcal{F}_{t}^{W}\right),
$$
where $x_{t,1}$ is the solution to our equation with interaction starting at time $t$ from the point $v$ and measure $\nu.$ Hence,
$$
E\left(\varphi^{\prime}(x(u, 1)) \mid \mathcal{F}_{t}^{W}\right)=\int_{\mathbb{R}} \varphi^{\prime}(v) \sigma_{t, 1}\left(\mu_{t}, x(u, t)\right)(d v),
$$
where $\sigma_{t, 1}(\nu, v)(\cdot)$ is the distribution of $x_{t, 1}(\nu, v).$
Due to the existence of the second stochastic derivative and the properties of the first and second derivatives for process $x $,
measure $\sigma_{t, 1}(\nu, v)$ has a differentiable density $p_{t, 1}(\nu, v, \cdot)$ with respect to Lebesgue measure.  Consequently,
$$
\begin{aligned}
E\left(\varphi^{\prime}(x(u, 1)) \mid \mathcal{F}_{t}^{W}\right) &=\int_{\mathbb{R}} \varphi^{\prime}(v) p_{t, 1}\left(\mu_{t}, x(u, t), v\right) d v \\
&=-\int_{\mathbb{R}} \varphi(v) \frac{\partial p_{t, 1}\left(\mu_{t}, x(u, t), v\right)}{\partial v} d v.
\end{aligned}
$$

Now consider the summand:

$$E(\varphi'(x(u,1))\int_t^1\alpha(u,s)\eta(u,s)ds\big | \mathcal{F}_t^W)=\int_t^1E(\varphi'(x(u,1))\eta(u,s)\alpha(u,s)\big | \mathcal{F}_t^W) ds.
$$

Similarly to the previous case,

\begin{eqnarray*}
&&E(\varphi'(x(u,1))\alpha(u,s)\eta(u,s)\big | \mathcal{F}_t^W)\\
&&\hskip0.7cm =\eta(u,s)\alpha(u,s)\int_{\mathbb{R}} \varphi'(v)\sigma_{s,1}(\mu_s,x(u,s),dv)\\
&&\hskip0.7cm =-\alpha(u,s)\eta(u,s)\int_{\mathbb{R}} \varphi(v)\frac{\partial}{\partial v} p_{s,1}(\mu_s,x(u,s),v)dv.
\end{eqnarray*}

Using Fubini theorem one can check that

$$E(\varphi(x(u,1))\int_t^1\alpha(u,s)\eta(u,s)ds\big | \mathcal{F}_t^W)=-\int_{\mathbb{R}}\varphi(v) E(\alpha(u,s)\eta(u,s)\frac{\partial}{\partial v} p_{s,1} (\mu_s,x(u,s),v)\big | \mathcal{F}_t^W) dv.
$$
Similarly to this case last  summand can be considered. Consequently the following theorem was proved.

\begin{thm}
\label{thm2.2}

Suppose that the coefficients  $a$ of equation with interaction

\begin{eqnarray}
&&d x(u, t)=\int_{\mathbb{R}}a\left(x(u, t),x(v,t)\right) \mu_{0}(dv)d t+ d W(t),\nonumber \\
&&x(u, 0)=u, u \in \mathbb{R},\nonumber \\
&&\mu_{t}=\mu_{0} \circ x^{-1}(u, t), t \geq 0,\label{7}
\end{eqnarray}
is bounded, together with two continuous partial derivatives by the positive constant $C$ such that

$$e^{3C}<1.$$

Then, the measure-valuated solution $\mu_t$ has Clark-Ocone's  representation which can be written as follows. For the test function $\varphi\in C_0^{\infty}$

$$\langle\varphi, {\mu}_{1}\rangle=E\langle\varphi, {\mu}_{1}\rangle+\int_0^1\int_{\mathbb{R}} \varphi(v) g_t(v)dW(t),$$
for some  random function $g$.
\end{thm}

Hence the integrant in Clark representation for measure-valued solution to equation with interaction takes values in the set of absolutely continuous measures.

\vskip 1 cm

\noindent
{\large\bf Jasmina \DJ or\dj evi\'c}\\[10pt]
Faculty of Science and Mathematics, University of Ni\v s, Vi\v segradska 33, 18000 Ni\v s, Serbia.\\
curraddr: The Faculty of Mathematics and Natural Sciences, University of Oslo, Blindern 0316 Oslo, Norway,\\
{\bf e-mai}l: jasmindj@math.uio.no, nina19@pmf.ni.ac.rs, djordjevichristina@gmail.com;
\vskip1cm

\noindent
{\large\bf Andrey Dorogovtsev} \\[10pt]
Institute of Mathematics National Academy of Sciences of Ukraine\\
{\bf e-mai}l: andrey.dorogovtsev@gmail.com.

\end{document}